\documentclass[a4paper]{amsart}
\usepackage{amssymb,amscd,amsthm,verbatim, latexsym}
\usepackage[english]{babel}
\newtheorem{theorem}{Theorem}[section]         
             
\newtheorem{corollary}[theorem]{Corollary}     
\newtheorem{proposition}[theorem]{Proposition} 
         
\newtheorem{remark}[theorem]{Remark}           
\newtheorem{definition}[theorem]{Definition}   
             
\newtheorem{example}[theorem]{Example}          
\numberwithin{equation}{section}

\newenvironment{namelist}[1]{%
\begin{list}{}
    {
      
      \settowidth{\labelwidth}{#1}
      \setlength{\leftmargin}{1.1\labelwidth}
    }
  }{%
\end{list}}
\newcommand{\nc}{\newcommand}
\nc{\R}{\ensuremath{\mathbb{R}}}
\nc{\K}{\ensuremath{\mathbb{K}}}
\nc{\Z}{\ensuremath{\mathbb{Z}}}
\nc{\N}{\ensuremath{\mathbb{N}}}
\nc{\Q}{\ensuremath{\mathbb{Q}}}
\nc{\Mo}{\ensuremath{\mathbb{M}}}
\nc{\C}{\ensuremath{\mathfrak {C}}}
\nc{\os}{\ensuremath{\mathfrak {S}}}
\nc{\A}{\ensuremath{\mathcal {A}}}
\nc{\E}{\ensuremath{\mathcal {E}}}
\nc{\D}{\ensuremath{\mathfrak {D}}}
\nc{\F}{\ensuremath{\mathfrak {F}}}
\nc{\B}{\ensuremath{\mathfrak {B}}}
\nc{\G}{\ensuremath{\mathfrak {G}}}
\nc{\rk}{\operatorname{rk}}
\nc{\IND}{\operatorname{IND}}
\nc{\NBC}{\operatorname{NBC}}
\nc{\Ker}{\operatorname{Ker}}
\nc{\OS}{\operatorname{OS}}
\nc{\OT}{\operatorname{OT}}
\nc{\M}{\ensuremath{\mathcal M}}
\nc{\I}{\ensuremath{\Im}}
\nc{\cl}{\ensuremath{c\ell}}
\nc{\cls}{\ensuremath{c\ell s}}
\nc{\clm}{\ensuremath{c\ell m}}
\nc{\sgn}{\operatorname{sgn}}
\nc{\id}{\operatorname{id}}

\begin{document}
\title[Diagonal bases]
{Diagonal bases in\\ Orlik-Solomon type algebras}
\thanks{ 
\emph{Secondary}:\, 05B35, 14F40.
}
\author{Raul Cordovil and David Forge}
\address{\newline 
Departamento de Matem\'atica,\newline 
Instituto Superior T\' ecnico \newline
 Av.~Rovisco Pais
 - 1049-001 Lisboa  - Portugal}
\email{cordovil@math.ist.utl.pt}
\thanks{The  first  author's research was 
supported in part by FCT (Portugal) through program POCTI and
 the project SAPIENS/36563/99.
The second  author's research  was  supported  by FCT trough the project
SAPIENS/36563/99.}
\address{{}\newline
Laboratoire de Recherche en Informatique\newline 
Batiment 490
Universite Paris Sud\newline 
91405 Orsay Cedex -
France
}
\email{forge@lri.fr}
\begin{abstract} To encode an important property of the ``no broken circuit bases"
of the Orlik-Solomon-Terao algebras, Andr\'as Szenes has introduced a particular type
of bases, the so called
 ``diagonal basis".  We
prove that this definition
 extends naturally to a large class of algebras, the so called
$\chi$-algebras.  Our definitions make also use of an ``iterative
residue formula" based on the matroidal operation of contraction. 
This formula can be seen as the combinatorial analogue of an iterative
residue formula introduced by Szenes.
As an application we deduce nice formulas to express
a pure  element 
in a diagonal basis. 
 \end{abstract}
\thanks{Typeset by \AmS -\LaTeX}
\maketitle 
\section{INTRODUCTION}
We denote by  $\mathcal{M}=\mathcal{M}([n])$
   a matroid 
 of rank $r$ on the ground set
$[n]:=\{1,2,\ldots,n\}$.
Let $V$ be a vector space of dimension $d$  over some field $\mathbb 
{K}$. A (central) arrangement (of 
hyperplanes) in $V,$
$\mathcal{A}_{\mathbb {K}}=\{H_{1},\ldots,H_{n}\},$ is a
   finite listed set of codi\-mension one vector
subspaces. 
Given an arrangement $\mathcal{A}_{\mathbb {K}}$\, we suppose always 
fixed a 
family of  linear forms 
$\big\{\theta_{{H}_i}\in V^{*}: H_i\in \mathcal{A}_{\mathbb
{K}},\,\mathrm{Ker}(\theta_{{H}_i})={H_i}
\big\},$ where\, $V^{*}$ denotes the dual space of\,
$V.$  
We denote by 
$L(\A_{\mathbb {K}})$  the {\em intersection lattice of
$\A_{\mathbb {K}}$}: i.e., the set of  intersections of 
hyperplanes in
$\mathcal{A}_{\mathbb {K}},$
 partially ordered by reverse inclusion.
There is a 
   matroid $\M(\A_{\mathbb {K}})$ on the ground set 
  $[n]$  determined by  
 $\mathcal{A}_{\K}$:
 a subset $D\subset [n]$\, 
 is  a \emph{dependent set} of  
$\mathcal{M}(\mathcal{A}_{\mathbb {K}})$ iff
 there are  scalars $\zeta_{i}\in \mathbb 
 {K},\, i\in D$, not all
 nulls,
 such that $\sum_{i\in D}\zeta_{i}\theta_{H_{i}}=0$. A \emph{circuit}  is a 
 minimal dependent set
with respect to inclusion. 

If $\K$ is an ordered field an additional 
structure is obtained: to every circuit $C,$ $\sum_{i\in C}\zeta_{i}\theta_{H_{i}}=0,$ we associate a
partition (determined up to a factor $\pm 1)$  $C^+=\{i\in C: \zeta_{i}>0\}, C^-=\{i\in
C:
\zeta_{i}<0\}.$ With this new structure $\M(\A_{\mathbb {K}})$ is said a 
\emph{(realizable) oriented matroid} and denoted by
$\boldsymbol{\M}(\A_{\mathbb {K}}).$ Set $\underline{\boldsymbol{\M}}(\A_{\mathbb {K}})=\M(\A_{\mathbb {K}}).$ Oriented
matroids on a ground set $[n]$, denoted $\boldsymbol{\mathcal{M}}([n]),$ are a very
natural mathematical concept and can be seen as the theory of
generalized hyperplane arrangements, see
\cite{MO}.
\par Set
$\mathfrak{M}(\mathcal{A}_{\mathbb{K}})=
V\setminus\bigcup_{H\in
 \mathcal{A}_{\mathbb{K}}} H.$ 
 The manifold \,$\mathfrak{M}(\mathcal{A}_{\mathbb{C}})$\, plays an important
role in the 
 Aomoto-Gelfand theory
 of multidimensional hypergeometric functions (see \cite{OT2} for a recent introduction
from the point of  view of arrangement theory).  Let ${K}$ be a commutative ring. In 
\cite{OS,OS1,OT}  the determination of the cohomology $K$-algebra
 $H^{*}\big(\mathfrak{M}(\mathcal{A}_{\mathbb{C}}); K\big)$\, from the matroid
 $\mathcal{M}(\mathcal{A}_{\mathbb {C}})$ is accomplished by first 
 defining the Orlik-Solomon $K$-algebra $\text{OS}(\mathcal{A}_{\mathbb{C}})$
 in terms of generators and 
 relators which depends only on the matroid $\mathcal{M}(\mathcal{A}_{\mathbb{C}})$ , and then by
showing that this
 algebra is isomorphic to
$H^{*}\big(\mathfrak{M}(\mathcal{A}_{\mathbb{C}}); K\big).$
The Orlik-Solomon 
algebras  have been then
intensively studied. Descriptions of developments from the early 1980's to the end of 1999, together with the contributions of many authors,
can be found in \cite{Falk1,Yuz}.

Aomoto suggested the study of the (graded) $\mathbb{K}$-vector space 
$\text{AO}
({\mathcal{A}}_{\mathbb{K}}),$ 
generated by the basis 
$\{Q({\mathcal{B}}_{\text{I}})^{-1}\},$ where $I$  is an independent
set of  
$\mathcal{M}
({\mathcal{A}}_{\mathbb{K}}),$   ${\mathcal{B}}_{\text {I}}:= \{H_i\in
\mathcal{A}_{\mathbb{K}}:i\in I\},$\, 
  and 
\,$Q
(
{\mathcal{B}}_{\text{I}}
)
=\prod_{i\in I}\theta_{H_i}$ denotes the corresponding defining polynomial. To 
answer to a conjecture of Aomoto,
Orlik and Terao have  introduced in \cite{OT1} a commutative 
$\mathbb{K}$-algebra, $\OT({\mathcal{A}}_{\mathbb{K}}),$ isomorphic to 
$\text{AO}(\mathcal{A}_{\mathbb {K}})$
 as a graded  $\mathbb {K}$-vector space in terms of the equations 
$\{\theta_{H}: H\in 
\mathcal{A}_{\mathbb{K}}\}.$\par
A  ``combinatorial
analogue'' of the   
  algebra of Orlik-Terao  was introduced in \cite{cor2}:
 to every oriented matroid $\boldsymbol{\mathcal{M}}$ was associated  a commutative
$\mathbb{Z}$-algebra, denoted by
$\mathbb{A}(\boldsymbol{\mathcal{M}}).$ \par
   Here we consider a large class of algebras, the so called $\chi$-algebras, that
contain the three just mentioned algebras: Orlik-Solomon, Orlik-Terao and the algebras
$\A_\chi(\boldsymbol{\mathcal{M}}),$ see \cite{FL} or Definition~\ref{FL1} below.
  Following Szenes \cite{Sz}, we define 
a particular type
of bases of $\A_\chi,$ the so called
 ``diagonal basis", see Definition
\ref{diagonal}. There is a natural example of these bases, the ``no circuit
basis".
We construct
the dual bases  of these bases, see Theorem~\ref{dbasis}.
 Our definitions make also use of an ``iterative
residue formula" based on the matroidal operation of contraction, see Equation
(\ref{resid}).  This formula can be seen as the combinatorial analogue of an iterative
residue formula introduced by Szenes, \cite{Sz}.
As applications we deduce nice formulas to express
a pure element
in a diagonal basis.
We prove also that the $\chi$-algebras verify a splitting short  exact
sequence, see Theorem~\ref{thm: There is}.
This theorem generalizes for the $\chi$-algebras previous similar theorems of
 \cite{cor2, OT}. 
\par
We use 
 \cite{W1, W2} as a general reference in matroid theory. We refer to 
 \cite{MO} and \cite{OT}
 for good sources
of the theory
of oriented matroids and arrangements of hyperplanes,  respectively.
%
\section{Diagonal bases}
 Let $\mathrm{IND}_\ell(\mathcal{M})\subset\binom{[n]}{\ell}$ be the family of 
independent sets  of cardinal $\ell$ of the matroid $\mathcal{M}$ and   set\,
$\mathrm{IND}(\mathcal{M}) =\bigcup_{\ell \in 
\mathbb{N}}\mathrm{IND}_\ell(\mathcal{M}).$  
We
denote by $\mathfrak{C}=\mathfrak{C}(\mathcal{M})$ 
  the set of 
circuits  of $\mathcal{M}$. For shortening of the notation
   the singleton set  $\{x\}$ is denoted  by  $x.$
When the smallest element $\alpha$ of a circuit $C,$ $|C|>1,$ is deleted, the
remaining  set,  $C\setminus \alpha,$\, is said to be a
\emph{broken circuit}. (Note that our definition is slightly different of the 
standard one. In
the standard definition
$C\setminus \alpha$ can be empty.) 
 A \emph{no broken circuit} set of a matroid $\mathcal{M}$ 
is an independent subset of $[n]$ which does not contain any broken circuit. 
Let
$\mathrm{NBC}_{\ell} (\mathcal{M})\subset\binom{[n]}{\ell}$ be the set of the 
no broken
circuit sets of cardinal $\ell$ of $\mathcal{M}.$  Set\, 
$\mathrm{NBC}(\mathcal{M})
=\bigcup_{\ell \in \mathbb{N}}\text{NBC}_\ell(\mathcal{M}).$  We denote by 
$L(\mathcal{M})$
 the lattice of  flats of
 $\mathcal{M}$.  \big(We remark that the 
 lattice map  
 $\phi: L(\mathcal{A}_{\mathbb {K}})\to L(\mathcal{M}
(\mathcal{A}_{\mathbb{K}}))$, determined by 
 the one-to-one correspondence 
 $\phi':H_{i}\longleftrightarrow \{i\}$, $i=1,\dotsc, n,$
 is a lattice isomorphism.\big{)}  For an   independent set $I,$  
let  $\cl(I)$ be the 
closure
of $I$ in $\mathcal{M}.$ \par
Fix a set $E=\{e_1,\ldots, e_n\}.$
Let  $\E=\K\oplus\E_1\oplus\cdots\oplus \E_n$ be the graded algebra over the field $\K$ generated by the
elements
$1, e_1,\dotsc,e_n$  and satisfying the relations
$e_i^2=0$ for all $e_i\in E$
and $e_j\cdot e_i=\beta_{i,j}e_i\cdot e_j$ with $\beta_{i,j}\in \K\setminus 0$ for all
$i<j.$ Both the  exterior algebra (take $\beta_{i,j}=-1$)
and the  commutative algebra with squares zero (take  $\beta_{i,j}=1$)
 are such algebras
and will be the only ones to be used in the examples.
 Let
$X^\sigma=(i_{\sigma(1)},i_{\sigma(2)},\,\ldots,\, i_{\sigma(m)}),$  $\sigma \in$
${\os}_m,$ denote the ordered set
$i_{\sigma(1)}\prec\cdots \prec i_{\sigma(m)}.$ 
When necessary we see the set $X=\{i_1,\dotsc, i_m\},$ as the ordered set
$X^{\mathrm{id}}.$ 
Set 
$X^\sigma\setminus
x:=(i_{\sigma(1)},\dotsc,\widehat{x},\dotsc,i_{\sigma(m)}
).$ If $Y^\beta=(j_{\beta(1)},\dotsc, j_{\beta(m')})$  and $X\cap Y=\emptyset,$ set
$X^\sigma*Y^\beta$ the concatenation $(i_{\sigma(1)},\dotsc,
i_{\sigma(m)},j_{\beta(1)},\dotsc, j_{\beta(m')}).$ In the sequel we will denote by
$e_{X}$ the   (pure) element $e_{i_1} e_{i_2} \cdots  e_{i_m}$
of\, $\mathcal E.$ Fix a mapping $\chi: 2^{[n]}\to \K.$  Let us also define $\chi$ for 
ordered sets by $\chi (X^\sigma )= \sgn(\sigma) \chi (X)$, where  $\sgn(\sigma)$ denotes 
the sign of the permutation $\sigma.$\par 
The $\chi$-{\it boundary}
of an element $e_X\in \E$  is   given by the equation
$$\partial e_X=\sum_{p=1}^{p=m}(-1)^p
\chi(X\setminus i_p)e_{X\setminus i_p}.$$
We extend $\partial$ to $\mathcal E$ by linearity. 
It is easy to see that for $\sigma \in \os_{|X|}$ we have $$\partial e_X=
\sgn(\sigma)\sum_{p=1}^{p=m}(-1)^p
\chi(X^\sigma \setminus i_{\sigma (p)})e_{X\setminus i_{\sigma (p)}}$$
and also for any $x\not \in X,$ $$\pm \partial e_{X\cup x}=
(-1)^{m+1}\chi (X)e_{X}+\sum_{p=1}^{p=m}(-1)^p
\chi(X\setminus i_{p}*x)e_{X\setminus i_{p}\cup x}.$$
Given an independent set $I,$ an element $a\in \cl(I)\setminus I$
is said 
\emph{active} in $I$ if $a$ is the minimal element of the unique circuit contained in
$I\cup a.$ We say that a subset $U\subset [n]$ is a \emph{unidependent}  of  $\M,$ if it
contains a unique circuit, denoted $C(U).$ 
 Note that
$U$ is unidependent iff
$\mathrm{rk}(U)=|U|-1.$ We say that an unidependent set $U$ is an \emph{inactive
unidependent} if $\mathrm{min} C(U)$ is the the minimal active element of 
$U \setminus\mathrm{min} C(U).$ 
Let us
remark that   $U$ is a unidependent of  $\M$ iff for some (or every) $x\in
U,$
$\mathrm{rk}(x)\not =0,$
$U\setminus x$ is a unidependent of $\M/x.$ 
\begin{definition}[\cite{FL}]\label{FL1}
{\em
Let ${\Im}_{\chi}(\M)$ be the (right) ideal of\, $\E$
generated by the\newline
$\chi$-boundaries $\{\partial e_C: C\in \C(\M), |C|>1\}$ and the set  
$\{e_i: \{i\}\in
\C(\M)\}.$  We say that
$\mathbb{A}_\chi(\M):=\E/{\Im}_{\chi}(\M)$
is a $\chi$-{\it algebra} if $\chi$ satisfies the
following two properties:
\begin{namelist}{xxxxxx}
\item[(UC1)]
$\chi(I)\not=0$ if and only if $I$ is independent.
\item[(UC2)]
For any two unidependents $U$  and $U'$ of $\M$ with $U'\subset U$ 
there is a  scalar $\varepsilon_{_{U,U'}}\in \K\setminus 0,$  such that
$\partial e_U=\varepsilon_{_{U,U'}}(\partial e_{U'})e_{U\setminus U'}.$
\end{namelist}
}
\end{definition}
\begin{remark}\label{rem1}
{\em
From (UC2) we conclude that ${\Im}_{\chi}(\M)$ has the basis
$$\{e_D: D ~\text{dependent of}~ \M\}\cup\{\partial e_U: U ~\text{inactive
unidependent of}~ \M\},$$
and that 
$\boldsymbol{nbc}:=\big\{[I]_{\mathbb{A}}:   I\in
\text{NBC}({\mathcal{M}})\big\}$ is a basis of 
 the vector space $\mathbb{A}=\mathbb{A}_\chi(\M).$  This 
fundamental property was first discovered for
the Orlik-Solomon algebras \cite{OT}, and
then also for other  classes of $\chi$-algebras, see \cite{cor2,
OT1}  and the following example for more details. Note also that this implies
that
$[X]_{\mathbb{A}}\not=0$ iff $X$  is an independent set of
${\mathcal{M}}.$
}
\end{remark}
\begin{example}[\cite{FL}]\label{chi}
{\em
Recall the three usual $\chi$-algebras. Let  $\E$ be the graded algebra over the field $\K$ generated by the
elements
$1, e_1,\dotsc,e_n$  and satisfying the relations
$e_i^2=0$ for all $e_i\in E$
and $e_j\cdot e_i=\beta_{i,j}e_i\cdot e_j$ where $\beta_{i,j}$ denotes a non null
scalar fixed for every pair 
$i<j.$ 
\begin{namelist}{xx}
\item[~~$\circ$] Let $\E$ be the  exterior algebra (taking $\beta_{i,j}=-1$).  Setting
$\chi(I^\sigma)=\text{sgn}(\sigma)$   for every independent set $I$ of a matroid $\M$ and
every permutation
$\sigma \in \os_{|I|},$ we obtain the Orlik-Solomon algebra,
$\OS(\M).$
\item[~~$\circ$]
Let $\A_\K=\{H_i: H_i=\Ker(\theta_i), i=1,2,\dotsc, n\}$ be an hyperplane arrangement
and
$\M(\A_\K)$  its associated  matroid.
For every flat $F:=\{f_1,\dotsc, f_k\}\subset [n]$ of $\M(\A_\K)$ we choose a basis $B_F$
of the vector subspace of $(\K^d)^*$ generated by 
$\{\theta_{f_1},\dotsc, \theta_{f_k}\}.$ By taking for $\mathcal E$ the free commutative
algebra with squares null  (taking $\beta_{i,j}=1$) and taking for any
$\{i_1,\dotsc,i_\ell\}=I\in
\mathrm{IND}_\ell,$\,\,
$\chi(I)=\det(\theta_{i_1},\dotsc,\theta_{i_\ell}),$ where the vectors are expressed in
the
  basis $B_{\cl(I)},$ we obtain the algebra $\OT(\A_\K),$
defined in
\cite{OT1}.
\item[~~$\circ$]
Let $\boldsymbol{\M}([n])$ be an oriented matroid. For every flat $F$ of $\boldsymbol{\M}([n]),$ we choose
(determined up to a factor $\pm 1$) a basis signature in the restriction of $\boldsymbol{\M}([n])$ to $F.$
We define a \emph{signature of the independents of an oriented matroid $\boldsymbol{\M}([n])$} as a mapping,
${\sgn}:\text{IND}(\boldsymbol{\M})\to \{\pm 1\},$ where $\sgn(I)$ is equal to the basis signature of $I$ in
the restriction of $\boldsymbol{\M}([n])$ to $\cl(I).$
By taking 
for $\mathcal E$ the free commutative algebra over the rational field $\Q$ with squares zero 
 (take  $\beta_{i,j}=1$) and taking $\chi(I)=\sgn(I)$ (resp.  $\chi(X)=0$)
 for every independent (resp. dependent) set of the matroid, we obtain the algebra
$\mathbb{A}(\boldsymbol{\M})\oplus_\Z \Q,$ where
$\mathbb{A}(\boldsymbol{\M})$ denotes the $\Z$-algebra defined in
\cite{cor2}.
\end{namelist}
}
\end{example}
For every $X\subset [n],$ we denote by $[X]_{\mathbb{A}}$ or
shortly by   $e_X$ when no confusion will result, the residue 
class in $\mathbb{A}_\chi({\mathcal{M}})$ 
determined by the  element $e_X.$  
Since ${\Im}_{\chi}(\M)$ is a homogeneous  ideal,   $\mathbb{A}_\chi(\M)$ inherits a grading from $\E.$ More precisely
we have
$\mathbb{A}_\chi({\mathcal{M}})=
\K\oplus\mathbb{A
}_1\oplus
\cdots 
\oplus 
\mathbb{A}_r,$ 
where $\mathbb{A}_\ell=\E_\ell/\E_{\ell}\cap{\Im}_{\chi}(\M)$ denotes the 
subspace  of $\mathbb{A}_\chi({\mathcal{M}})$ generated by the elements 
$\big\{[I]_{\mathbb{A}}: I\in 
\text{IND}_\ell({\mathcal{M}})\big\}.$
Set $\boldsymbol{nbc}_{\,\ell}:=\big\{[I]_{\mathbb{A}}:  
I\in \text{NBC}_{\,\ell}({\mathcal{M}})\big\}$ and
$\boldsymbol{nbc}:=\bigcup_{\ell=0}\boldsymbol{nbc}_{\,\ell}.$
From Remark \ref{rem1} we conclude that
$\boldsymbol{nbc}_\ell$ is
a   basis of the vector space $\mathbb{A}_\ell.$  
\begin{proposition}\label{con}
Let $\mathbb{A}_\chi({\M})$ be a $\chi$-algebra. For any non loop element $x$ of
$\M([n]),$ we define the two maps:
\begin{equation}\label{eq1}
\chi_{\M \setminus x} : 2^{[n]\setminus x}\to \K\hspace{5mm}  \text{by}\hspace{5mm}
\chi_{\M \setminus x}(I)=\chi(I)\hspace{5mm}  \text{and}
\end{equation} 
\begin{equation}\label{eq2}
\chi_{\M /x} : 2^{[n]\setminus x}\to \K \hspace{5mm}  \text{by}\hspace{5mm}
\chi_{{\M /x}}(I)=\chi(I*x).
\end{equation} 
 Then \,\,
$\mathbb{A}_{\chi_{\M /x}}({\M/x})$\,\, and\,\, 
$\mathbb{A}_{\chi_{\M \setminus x}}({\M\setminus x})$\,\, are\,\,  $\chi$-algebras.
\end{proposition}
\vspace{2mm}
\noindent{\em Proof}.\,
The deletion case being trivial, we will just prove the contraction case.
We have to show that $\chi_{\M /x}$  verifies properties (UC1) and (UC2).
The first property is verified since a set $I$ is independent in  $\M /x$ iff
$I\cup x$ is independent in $\M.$ To see that the second property is also verified, let
$U$ and $U'$ be two unidependents of $\M /x$ (iff $U\cup x$ and $U'\cup x$ 
are two unidependents of $\M$). We
know that $\partial e_{U\cup x}=\varepsilon_{_{{U\cup x},{U'\cup x}}}
(\partial e_{{U'\cup x}})e_{U\setminus U'}.$
We note $\partial^\prime$ the boundary defined by  $\chi_{\M /x}$ and so we will show
that there is a scalar $\varepsilon_{_{U,U'}}$ such that 
$\partial^\prime e_U=\varepsilon_{_{U,U'}}(\partial^\prime e_{U'})e_{U\setminus U'}.$
Let $X,X'\subset [n]$ two disjoint subsets then $e_Xe_{X'}=\beta_{_{X,X'}e_{X\cup
X'}},$ where $\beta_{_{X,X'}}=\prod _{e_i\in X,e_j\in X',\, i>j}\beta_{i,j}.$ We have
with $U=(i_1,\dotsc, i_{m})$ and $U'=(j_1,\dotsc, j_{k})$:
 $$\pm\partial e_{U \cup x}= \sum_{p=1}^{p=m}(-1)^p
\chi(U\setminus i_p*x)e_{U\cup x\setminus i_p}+(-1)^{m+1}\chi(U)e_U,$$
$$\partial^\prime e_{U}= \sum_{p=1}^{p=m}(-1)^p
\chi(U\setminus i_p*x)e_{U\setminus i_p},$$
\begin{align*}
\pm(\partial e_{U' \cup x})e_{U\setminus U'}
=\sum_{p=1}^{p=k}(-1)^p
\chi(U'\setminus j_p* x)\beta_{_{U'\cup x\setminus j_p, U\setminus U'}}
e_{U\cup x\setminus j_p}+\\
+(-1)^{k+1}\chi(U')\beta_{_{U',U\setminus U'}}e_{U},
\end{align*}
 $$(\partial^\prime e_{U' })e_{U\setminus U'}=\sum_{p=1}^{p=k}(-1)^p
\chi(U'\setminus j_p*x)\beta_{_{U'\setminus j_p,U\setminus U'}}e_{U\setminus j_p}.$$
After remarking that
$\beta_{_{U'\cup x\setminus j_p,U\setminus U'}}
\beta_{_{U'\setminus j_p,U\setminus U'}}^{-1}=\beta_{_{x,U\setminus U'}}$  does not
depend on $j_p,$ we can deduce that 
$\partial' e_U=\varepsilon_{_{U,U'}}(\partial' e_{U'})e_{U\setminus U'}$
with $\varepsilon_{_{{U\cup x},{U'\cup x}}}=\pm\varepsilon_{_{U,U'}}\beta_{x,U\setminus
U'}.$\qed
\vspace{2mm}
\begin{proposition}\label{l} For every    non loop element $x$ of
$\M([n]),$
there is a unique monomorphism of\, vector spaces,\,
$\mathfrak{i}_x:  \mathbb{A}(\M\setminus x)\to
\mathbb{A}(\M),$ such that such that, for every
$I\in
\mathrm{IND}(\M\setminus x),$ we have
$\mathfrak{i}_x(e_I)=e_I.$  
  \end{proposition}
\begin{proof}
 By a 
reordering  of the elements of the matroid
$\M$ we can suppose
that
     $x=n.$
It is clear that 
$$\text{NBC}(\M\setminus x)=\big\{X:
     X\subset [n-1] 
~~\text{and}~~X\in \text{NBC}(\M)\big\},$$
so the proposition is a consequence of Equation~(\ref{eq1}).
\end{proof}
\begin{proposition} For every non loop element $x$ of
$\M([n]),$ 
there is a unique epimorphism of vector spaces, 
$\boldsymbol{\mathfrak{p}}_x: 
\mathbb{A}({\mathcal{M}})\to
\mathbb{A}({\mathcal{M}}/x),$
 such that, for every
$e_I,$ $I\in
\mathrm{IND}({\mathcal{M}}),$ we have
\begin{equation}\label{p_x}
\boldsymbol{\mathfrak{p}}_x(e_I):=
\begin{cases}
\vspace{2mm}
e_{I\setminus x} & 
   \mbox{\em if} \hspace{5mm} x\in I,\\
\frac{\chi ( I\setminus y , x)}{\chi (I\setminus y,y)}\,e_{I\setminus y} & 
   \mbox{\em if there is }  y\in I 
 \hspace{2mm}\mbox{\em  parallel 
to $x,$}\\
0 & \mbox{\em otherwise}. 
\end{cases}
\end{equation}
\end{proposition}
\noindent{\em Proof.}\, From Remark~\ref{rem1},  it is enough to
prove that $\mathfrak{p}_x({\partial}e_{U})=0,$ 
for all unidependent $U=(i_1,\dotsc,i_{m}).$
We recall that if $x\in U$ then
${{U}}\setminus x$ is a unidependent set of
${\mathcal{M}}/x.$ 
There are only the  following four cases:
\begin{namelist}{xx}
\item[~~$\circ$] If $U$ contains $x$ but no $y$ parallel to $x$ then:
\begin{align*}
\pm\mathfrak{p}_x({\partial}e_{U})=\mathfrak{p}_x(
(-1)^{m}\chi(U\setminus x)e_{U\setminus x}+
\sum_{i_p\in U\setminus x}(-1)^p
\chi(U\setminus \{i_p,x\}*x)e_{U\setminus i_p}))\\
=\sum_{i_p\in U\setminus x}(-1)^p
\chi(U\setminus \{i_p,x\}*x)e_{U\setminus \{i_p,x\}}=0\hspace{30mm}
\end{align*}
from Proposition \ref{con}.
\item[~~$\circ$] If $U$ does not 
contain $x$ but a $y$ parallel to $x$ then:
\begin{align*}
\pm\mathfrak{p}_x({\partial}e_{U})=\mathfrak{p}_x\big(
(-1)^{m}\chi(U\setminus y)e_{U\setminus y}+
\sum_{i_p\in U\setminus y}(-1)^p
\chi(U\setminus \{i_p,y\}*y)e_{U\setminus i_p}\big)\\
=
\sum_{i_p\in U\setminus y}(-1)^p
\chi(U\setminus \{i_p,y\}*y)\frac{\chi(U\setminus \{i_p,x\}*x)}
{\chi(U\setminus \{i_p,y\}*y)}e_{U\setminus \{i_p,y\}}=
0\hspace{1mm}
\end{align*}
like previously since $U\setminus y$ is again a unidependent of $\M/x.$
\item[~~$\circ$] If $U$   contains $x$ and a $y$ parallel
  to $x$ then:
\begin{align*}
\pm\mathfrak{p}_x({\partial}e_{U})=\mathfrak{p}_x(
\chi(U\setminus\{x,y\}*y)e_{U\setminus x}-\chi(U\setminus\{x,y\}*x)e_{U\setminus
y})\\
=
\chi(U\setminus \{x,y\}*y)\frac{\chi(U\setminus \{x,y\}*x)}
{\chi(U\setminus \{x,y\}*y)}e_{U\setminus 
\{x,y\}}-\chi(U\setminus\{x,y\}*x)e_{U\setminus  
\{x,y\}}=
0.
\end{align*}
 \item[~~$\circ$] If $U$ does not 
contain $x$ nor a $y$ parallel to $x$ then:
\begin{equation*}
\hspace{25mm}\mathfrak{p}_x({\partial}e_{U})=\mathfrak{p}_x\big(
\sum_{i_p\in U}(-1)^p
\chi(U\setminus i_p)e_{U\setminus i_p}\big)=0.
\hspace{25mm}\qed
\end{equation*}
\end{namelist}
 \begin{theorem}\label{thm: There is}
For every  element $x$ of a simple $\M([n]),$    there is a splitting short  exact
sequence of vector spaces
\begin{equation}\label{exact}
    0\to \mathbb{A}(\M\setminus
x)\stackrel{\mathfrak{i}_x}\longrightarrow 
    \mathbb{A}(\M)\stackrel{\mathfrak{p}_x}\longrightarrow 
    \mathbb{A}(\M/x)\to 0.
    \end{equation}
\end{theorem}
\begin{proof}
From the definitions we know that $\mathfrak{p}_x\circ\mathfrak{i}_x$, is the null map
 so $\mathrm{Im}(\mathfrak{i}_x)\subset \mathrm{Ker}(\mathfrak{p}_x).$ 
    We will prove the equality  $\dim(\mathrm{Ker}(\mathfrak{p}_n))=\dim
(\mathrm{Im}(\mathfrak{i}_n)).$ By a reordering  of the elements of $[n]$
 we can suppose
that
     $x=n.$
The minimal 
 broken circuits of 
$\M/n$ are the minimal sets $X$ such that either $X$
or $X\cup \{n\}$ is a broken circuit of $\M$ (see the Proposition 3.2.e
of \cite{Bry}).  Then
\begin{equation*}
\text{NBC}(\M/n)=\big\{X: X\subset 
 [n-1]~~\mbox{and}~~X\cup\{n\}\in \text{NBC}(\M)
\big\}\,\,\,\,\,\,\mbox{and}
\end{equation*}
\begin{equation}\label{NBC}
\text{NBC}(\M)=
\text{NBC}(\M\setminus n)\biguplus
\big\{I\cup
n: I\in
\text{NBC}(\M/n)\big\}.
\end{equation}
So $\dim(\mathrm{Ker}(\mathfrak{p}_n))=\dim
(\mathrm{Im}(\mathfrak{i}_n)).$ There is a  
 morphism of modules 
$$\mathfrak{p}^{-1}_n: \mathbb{A}(\M/n)\to
\mathbb{A},~~~\text{where}~ ~~~\mathfrak{p}^{-1}_n([I]_{\mathbb{A}(\M/n)}):=[I\,\cup
\,n]_{\mathbb{A}}, \forall I\in
\mathrm{NBC}(\M/n).$$
It is clear that $\mathfrak{p}_n\circ \mathfrak{p}^{-1}_{n}$ is the
identity map. From Equation~(\ref{NBC}) we conclude that the exact
sequence~(\ref{exact}) splits.
\end{proof}
Similarly to  \cite{Sz} (see also \cite{BV}), we now construct, making use  
of iterated contractions,  the dual basis 
 $\boldsymbol{nbc}_{\,\ell}^*=(b^*_i)$ of the basis  $\boldsymbol{nbc}_{\,\ell}=(b_j).$ 
 More precisely
$\boldsymbol{nbc}^*_{\,\ell}$ is the basis of $\mathbb{A}_\ell^*$ 
the vector space of the linear forms such that 
$b^*_i(b_j)=\delta_{ij}$ (the Kronecker delta).

We associate to  the ordered independent set 
$I^\sigma:=(i_{\sigma(1)},\dotsc,  i_{\sigma(p)})$ of\,
${\mathcal{M}}$ the linear form on 
$\mathbb{A}_\ell,$ 
$\boldsymbol{\mathfrak{p}}_{I^\sigma}:
\mathbb{A}_\ell\to  \K,$
\begin{equation}\label{resid}
\boldsymbol{\mathfrak{p}}_{I^\sigma}:=
\boldsymbol{\mathfrak{p}}_{e_{i_{
\sigma(1)}}}\!\circ\boldsymbol{\mathfrak{p}}_{e_{i_{\sigma(2)}}}\!\circ
\cdots\,\circ
\boldsymbol{\mathfrak{p}}_{e_{i_{\sigma(p)}}}.
\end{equation}
   We call
$\boldsymbol{\mathfrak{p}}_{I^\sigma}$  the  \emph{iterated residue} with 
respect to the ordered independent
set $I^\sigma.$ (It is clear that the map $\boldsymbol{\mathfrak{p}}_{I^\sigma}$ depends on the 
order  chosen on $I^\sigma$ and not only on
the  underlying set $ I.)$ We
associate to $I^\sigma$ the flag of flats of
${\mathcal{M}},$
\begin{equation*}
\textbf{Flag}(I^\sigma):=\cl\big(\{i_{\sigma(p)}\}\big) \subsetneq 
\cl\big(\{i_{\sigma(p)}, i_{\sigma(p-1)}\}\big)
\subsetneq 
\dotsb   \subsetneq  \cl(I).
\end{equation*}
\begin{proposition}\label{resi}
Let $J\in 
\mathrm{IND}_\ell({\mathcal{M}})$ then we have 
$\boldsymbol{\mathfrak{p}}_{I^\sigma}(e_{J})\not=0$ iff there is a unique 
permutation $\tau \in \os_\ell$  such that
$\textbf{Flag}\,(J^\tau) =
\textbf{Flag}\,(I^\sigma).$  And in this case we have 
$\boldsymbol{\mathfrak{p}}_{I^\sigma}(e_{J})=
\chi (I^\sigma )/\chi (J^\tau).$ In particular we have
$\boldsymbol{\mathfrak{p}}_{I^\sigma}(e_{I})=1$ for any independent set $I$ 
and any permutation $\sigma.$
\end{proposition}
\begin{proof}
The first equivalence is very easy to prove in both direction.
To obtain the expression of $\boldsymbol{\mathfrak{p}}_{I^\sigma}(e_{J})$
we just need to iterate $\ell$ times the residue. This gives:
\begin{align*}
\boldsymbol{\mathfrak{p}}_{I^\sigma}(e_{J})=
\frac{\chi ( J\setminus j_{\tau (\ell)} *i_{\sigma(\ell)})}
{\chi ( J\setminus j_{\tau (\ell)}*j_{\tau (\ell)} )}
\times \frac{\chi ( J\setminus \{j_{\tau (\ell)}, j_{\tau (\ell
-1)}\}*  i_{\sigma (\ell -1)}*i_{\sigma(\ell)})}
{\chi ( J\setminus \{j_{\tau (\ell)}, j_{\tau (\ell -1)}\}*
j_{\tau (\ell-1)}*i_{\sigma(\ell )})} \times \dotsb \\
\dotsb \times
\frac{\chi (I^\sigma)}{\chi (j_{\tau (1)}*I^\sigma \setminus 
i_{\sigma(1)})}.
\end{align*}
After simplification we obtain the announced formula. And finally the
last result comes from the fact that if $I=J$ then clearly $\tau=\sigma.$
\end{proof}
\begin{remark}\label{rem}
{\em
The fact that $\boldsymbol {\mathfrak{p}}_{I^\sigma}(e_{J})$ is null 
depends on the permutation $\sigma$. For example, for any
simple matroid of rank 2  we have  
$\boldsymbol {\mathfrak{p}}_{13}(e_{12}) =0$ and 
$\boldsymbol{\mathfrak{p}}_{31}(e_{12})\not =0$.  But if $\boldsymbol
{\mathfrak{p}}_{I^\sigma}(e_{J})
\not=0$ then its value does not depend on $\sigma.$ We mean by this that if there 
are two permutations $\sigma $ and $\sigma '$ such that 
$\boldsymbol{\mathfrak{p}}_{I^\sigma}(e_{J})\not =0$
and $\boldsymbol{\mathfrak{p}}_{I^{\sigma '}}(e_{J})\not =0$ then 
$\boldsymbol{\mathfrak{p}}_{I^\sigma}(e_{J})=\boldsymbol{\mathfrak{p}}_{I^{\sigma
'}}(e_{J}).$
}
\end{remark}
\begin{definition}[\cite{Sz}]\label{diagonal}
 {\em
We
say that the subset $\mathbb{I}_\ell\subset \big\{[I]_{\mathbb{A}}: I\in 
\text{IND}_\ell({\mathcal{M}})\}$ is
 a
\emph{diagonal basis} of   $\mathbb{A}_\ell$ if and only if the following three
conditions hold:
\begin{namelist}{xxxxxx}
\item[$(\ref{diagonal}.1)$] For every $[I]_{\mathbb{A}}\in \mathbb{I}_\ell$ there is a fixed
permutation of the set $I$ denoted $\sigma_I\in \os_\ell;$
 \item[$(\ref{diagonal}.2)$] $\big|\mathbb{I}_\ell| 
\geq\text{dim}(\mathbb{A}_\ell);$
\item[$(\ref{diagonal}.3)$]  For every $[I]_{\mathbb{A}}, [J]_{\mathbb{A}}\in
\mathbb{I}_\ell$ and every permutation
$\tau
\in \os_\ell,$ the equality
$\textbf{Flag}\,(J^\tau) =
\textbf{Flag}\,(I^{\sigma_I})$ implies $J= I.$
\end{namelist}
}
\end{definition}
\begin{theorem}\label{dbasis} Suppose that $\mathbb{I}_\ell$ is a diagonal 
basis  of   $\mathbb{A}_\ell.$
Then $\mathbb{I}_\ell$ is a basis
of  $\mathbb{A}_\ell$ and
\,$\mathbb{I}_\ell^*:=\{\boldsymbol{\mathfrak{p}}_{I^{\sigma_I}}: 
[I]_\mathbb{A}\in \mathbb{I}_\ell\}$ is the dual 
 basis  of\, $\mathbb{I}_\ell.$ 
\end{theorem}
\begin{proof} 
Pick two elements $[I]_\mathbb{A}, [J]_\mathbb{A}
\in\mathbb{I}_\ell.$ Note that
$\mathfrak{p}_{I^{\sigma_I}}(e_{J})=\delta_{IJ}$ (the Kronecker delta), from Condition~$(\ref{diagonal}.2)$ and
Proposition~\ref{resi}.
 The elements of $\mathbb{I}_\ell$ are linearly independent: suppose that
$[J]=\sum\zeta_j [I_j],\,  \zeta_j\in \K\setminus 0;$ then
$1=\boldsymbol{\mathfrak{p}}_{J^{\sigma_J}}([J])=
\boldsymbol{\mathfrak{p}}_{J^{\sigma_J}}\big(\sum\zeta_j
[I_j]\big)=0,$ a contradiction. It is clear also that $\mathbb{I}_\ell^*$ 
is the dual 
 basis  of\, $\mathbb{I}_\ell.$ 
\end{proof}
The following result gives an interesting explanation of results of \cite{cor} and \cite{cor2}.
\begin{corollary}\label{that}
$\boldsymbol{nbc}_\ell({\mathcal{M}})$  is a diagonal basis of
$\mathbb{A}_\ell$ where $\sigma_I$ is the identity  for every $[I]_{\mathbb{A}}
\in
\boldsymbol{nbc}_\ell({\mathcal{M}}).$
  For a given $[J]_{\mathbb{A}}\in \mathbb{A}_\ell,$
suppose
that
\begin{namelist}{xxxxxxxx}
\vspace{1mm}
\item[$(\ref{that}.2)$]$[J]_{\mathbb{A}}
=\sum
\xi (I,J) [{I}]_{\mathbb{A}},$ where
\,\,$[I]_{\mathbb{A}}\in\boldsymbol{nbc}_\ell({\mathcal{M}})$\, and 
\,\, 
$\xi (I,J) \in \K.$
\vspace{1mm}
\end{namelist}
Then  are equivalent:
\begin{namelist}{xxxxx}
\item[~~$\circ$]$\xi (I,J)\not =0,$
\item[~~$\circ$]$\textbf{Flag}\,(I) =
\textbf{Flag}\,({J}^{\tau})$ for some permutation $\tau.$
\end{namelist} 
If $\xi (I,J)\not =0$ we have $\xi (I,J)=\frac{\chi (I)}{\chi (J^\tau)}.$ In particular
if $\mathbb{A}$  is the Orlik-Solomon algebra then
$\xi(I,J)=\sgn(\tau).$
\end{corollary}
\begin{proof}
By hypothesis $(\ref{diagonal}.1)$ and $(\ref{diagonal}.2)$ are true.
 We claim that 
$\boldsymbol{nbc}_{\,\ell}({\mathcal{M}})$ verifies
$(\ref{diagonal}.3).$ Suppose for a contradiction that $J\not = I,$
$[J]_\mathbb{A}, [I]_\mathbb{A}\in 
\boldsymbol{nbc}_{\,\ell}({\mathcal{M}})$ and there is $\tau \in 
\os_\ell,$ such that $\textbf{Flag}\,(J^\tau) =
\textbf{Flag}\,(I).$ Set $I=(i_1,\dotsc,  i_\ell)$ and  
$J=(j_{\tau(1)},\dotsc,  j_{\tau({\ell})}),$  and suppose that $j_{\tau(m+1)}=i_{m+1}, 
\dotsc, j_{\tau(\ell)}=i_\ell$ and $i_m\not =j_{\tau({m})}.$ Then there is a circuit
$C$ of
$\M$ such that $$i_m, j_{\tau(m)}\in 
C\subset \{i_m, j_{\tau(m)}, i_{m+1}, i_{m+2}, \dotsc
,i_\ell\}.$$
If $j_{\tau(m)}<i_m$ [resp. $i_m<j_{\tau(m)}$] we conclude that $I\not \in 
\text{NBC}_\ell({\mathcal{M}})$
[resp. $J\not \in \text{NBC}_\ell({\mathcal{M}})$] a 
contradiction. So $\boldsymbol{nbc}_{\,\ell}({\mathcal{M}})$  is a diagonal basis of
$\mathbb{A}_\ell.$\par
From Theorem~\ref{dbasis} we 
conclude that   
$\boldsymbol{nbc}_\ell^*:=\big\{\boldsymbol{\mathfrak{p}}_{I}: 
[I]_\mathbb{A}\in
\boldsymbol{nbc}\}$ is the dual 
 basis  of\, $\boldsymbol{nbc}.$ Suppose now that 
$[J]_{\mathbb{A}}=\sum{\xi}_{I}[{I}]_{\mathbb{A}},$ where
$[I]_{\mathbb{A}}\in\boldsymbol{nbc}_\ell({\mathcal{M}})$ and 
${\xi}_I\in
k.$ Then ${\xi}_{I}=\boldsymbol{\mathfrak{p}}_{I}(e_{J})$
and the remaining follows from Proposition~\ref{resi}.
\end{proof}
Making full use of the matroidal notion of iterated residue, see Equation~(\ref{resid}),
we are able to prove the following result very close to Proposition~ 2.1 of
\cite{Sz2}.
\begin{proposition}\label{machain}
Consider the set  of  vectors  $\mathcal{V}:=\{v_1,\dotsc,  v_k\}$ 
in the plane $x_d=1$ of\, $\K^d.$ Set $\A_\K:=\{H_i: H_i=\Ker(v_i)\subset (\K^d)^*,\,
i=1,\dotsc, k\}$  and let\, $\OT(\A_\K)$ be its Orlik-Terao 
corresponding
algebra.  Fix a diagonal basis  
$\mathbb{I}_\ell\subset \{[I]_{\mathbb{A}}: I\in \mathrm{IND}_\ell(\M)\}
$
of  $\mathbb{A}_\ell$ and let\,
$\mathbb{I}^*_\ell=\{{\mathfrak{p}}_{I^{\sigma_I}}: [I]_{\mathbb{A}}\in \mathbb{I}_\ell\}
$
be the corresponding dual  basis. 
Then, for any\,  $e_J\in \mathbb{A}_\ell \setminus 0,$ 
we have 
$$\sum _{I\in \mathbb{I}_\ell} \boldsymbol {\mathfrak{p}}_{I^{\sigma_I}}(e_J)=\sum _{I\in
\mathbb{I}_\ell}
\big\langle {\mathfrak{p}}_{I^{\sigma_I}},e_J\big\rangle =1.$$
\end{proposition}
\begin{proof}
We have for any $\ell+1$-subset of $\mathcal{V},$
$\sum_{p=1}^{p=\ell+1}(-1)^p
\chi(U\setminus i_p)=0.$ (This is the development of a determinant with two lines of 
1.)
For any rank $\ell$ unidependent $U=\{i_1,\dotsc,  i_{\ell+1}\}$  of the matroid
$\mathcal{M} (\A_\K),$ we have 
$$\partial e_U=\sum_{p=1}^{p=\ell+1}(-1)^p\chi(U\setminus i_p)e_{U\setminus i_p}.$$
Since  the sum of the coefficients in these relations is 0 and that these relations are
generating, see Remark~\ref{rem1}, we can deduce that the sum of the coefficients in any
relation in
$\OT(\A_\K)$ is also equal to 0 which concludes the proof.
\end{proof}
\begin{example}
{\em Consider  the  6 points $p_1,\dotsc,  p_6$ in the
affine plane $z=1$ of
$\R^3,$  whose coordinates are indicated in 
Figure~1. Set
$v_i:=\overrightarrow{(0,p_i)},\ i=1,\dotsc, 6.$ And let
$\A$ be
 the corresponding arrangement  of
$(\R^3)^*,$ $\A:=\{H_i=\Ker(v_i),\, i=1,\dotsc, 6\}.$
Let $\M(\A)$ [resp. $\boldsymbol{\M}(\A)$] be the corresponding rank three
[resp. oriented]  matroid.
\begin{center}
\begin{picture}(-76,120)(100,10)
\put(15,15){$p_1$}\put(15,110){$p_3$}\put(15,63){$p_2$}\put(63,15){$p_4$}\put(112,15){$p_5$}
\put(34,34)
{$p_6$}
\put(10,10) {\circle*{6}}
\put(60,10) {\circle*{6}}
\put(10,60) {\circle*{6}}
\put(10,110) {\circle*{6}}
\put(44,43) {\circle*{6}}
\put(110,10) {\circle*{6}}
\put(-25,10){(0,0,1)}\put(-28,60){$(0,\frac {1}{2},1)$}\put(-25,110){(0,1,1)}
\put(43,-2){$(\frac {1}{2},0,1)$}\put(97,-2){(1,0,1)}\put(47,49){$(\frac {1}{3},\frac
{1}{3},1)$}
\put(10,10) {\line(2,0){100}}
\put(10,10) {\line(0,1){100}}
\put(10,60) {\line(2,-1){100}}
\put(10,110) {\line(1,-2){50}}
\end{picture}
\end{center}
\vspace{8mm}
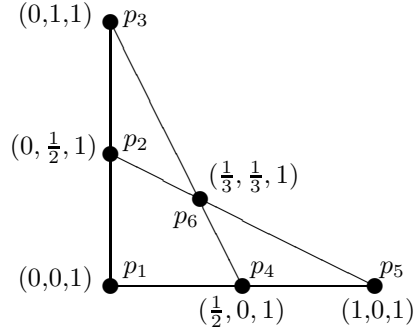
\centerline{\textbf{Figure 1}}
\vspace{3mm}
Let $\mathbb{A}_\chi$ be a  $\chi$-algebra on $\M(\A).$
We know that $$\boldsymbol{nbc}_{\,3}=\{e_{124},e_{125},e_{126},
e_{134},e_{135},e_{136}\}$$ together with $\sigma _{124}=\sigma _{125}=
\sigma _{134}=\sigma _{135}=\sigma _{136}=\sigma _{156}=
\text{id}$ is a diagonal basis of $\mathbb{A}_3,$ from Corollary~\ref{that}.
Directly from  the
Definition~\ref{diagonal} we see that  $\B_3=\{ e_{124},e_{125},e_{134},
e_{135},e_{136},e_{156}\}$ with
$\sigma _{124}=\sigma _{134}=\sigma _{135}=\sigma _{136}=\sigma _{156}=
\text{id}$ and $\sigma _{125}=(132)$ is also a diagonal basis of 
$\mathbb{A}_3.$ We will look at expression 
on the  basis $\boldsymbol{nbc}_3$ (resp. $\B_3$) of the vector space 
$\mathbb{A}_3,$
of some
elements of the type $e_B,$ $B$ basis of $\M(\A),$   for the three 
$\chi$-algebras  of  Example~\ref{chi}. Especially, one can verify   as
 stated in Remark~\ref{rem}\, that
 $\boldsymbol{\mathfrak{p}}_{125^{\id}}(e_{235})=
\boldsymbol{\mathfrak{p}}_{125^{(132)}}(e_{235}).$ Let also point out that for
 the Orlik-Terao algebra, we have $\sum _{I\in
 \B}\boldsymbol{\mathfrak{p}}_{I^{\sigma}}(e_{J})=1$
as proved in Proposition~\ref{machain}.
\begin{namelist}{xx}
\item[~~$\circ$]Consider the  basis $\boldsymbol{nbc}_3$ of the $\mathbb{K}$-vector space
$\mathbb{A}_3.$ So
 we have:
$$e_{235}=\sgn(325)e_{125} + \sgn(235)e_{135}=
-e_{125} + e_{135}~~ \text{in}~~\OS (\M(\A)),$$
$$e_{235}=\frac {\det (125)}{\det (325)}e_{125}+ 
\frac {\det (135)}{\det (235)}e_{135}=-e_{125}+2e_{135}~~~ \text{in}~~ 
\OT(\A),$$
$$e_{235}=\chi (125)\chi (325)e_{125}+ 
\chi (135)\chi (235)e_{135}=-e_{125}+e_{135}~~~ \text{in}~~~
\mathbb{A}(\boldsymbol{\M}(\A)).$$ 
$$e_{156}=\text{sgn}(165)e_{125}+ \text{sgn}(156)e_{126}=
-e_{125}+e_{126}~~~\text{in}~~~ \OS (\M(\A)),$$
$$e_{156}=\frac {\det (125)}{\det (165)}e_{125}+ 
\frac {\det (126)}{\det (156)}e_{126}=\frac {3}{2}e_{125}-\frac {1}{2}e_{126}
~~~\text{in}~~~\OT(\A),$$
$$e_{156}=\chi (125)\chi (165)e_{125}+ 
\chi (126)\chi (156)e_{126}=e_{125}-e_{126}~~~\text{in}~~~
\mathbb{A}(\boldsymbol{\M}(\A)).$$
\vspace{-2mm}
\item[~~$\circ$]Consider now the  basis $\B_3$ of the $\mathbb{K}$-vector space
$\mathbb{A}_3.$ So we have:
$$e_{235}=\text{sgn}(152)\text{sgn}(352)e_{125}+ \text{sgn}(235)e_{135}=
-e_{125}+e_{135}~~~ \text{in}~~ \OS (\M(\A)),$$
$$e_{235}=\frac {\det (152)}{\det (352)}e_{125}+ 
\frac {\det (135)}{\det (235)}e_{135}=-e_{125}+2e_{135}~~~ \text{in}~~\OT(\A),$$
$$e_{235}=\chi (152)\chi (352)e_{125}+ 
\chi (135)\chi (235)e_{135}=-e_{125}+e_{135}~~~
\text{in}~~\mathbb{A}(\boldsymbol{\M}(\A)).$$ 
$$e_{126}=\text{sgn}(162)\text{sgn}(152)e_{125}+ \text{sgn}(126)e_{156}=
e_{125}+e_{156}~~~ \text{in}~~\OS (\M(\A)),$$
$$e_{126}=\frac {\det (152)}{\det (162)}e_{125}+ 
\frac {\det (156)}{\det (126)}e_{156}=3e_{125}-2e_{156}~~ \text{in}~~\OT(\A),$$
$$e_{126}=\chi (152)\chi (162)e_{125}+ 
\chi (156)\chi (126)e_{156}=e_{125}-e_{156}~~
\text{in}~~ \mathbb{A}(\boldsymbol{\M}(\A)).$$ 
\end{namelist}
}
\end{example}

\end{document}